\documentclass[11pt,reqno]{amsart}
\usepackage{amsfonts,amssymb,amsmath}
\usepackage{color}

\DeclareMathOperator{\diag}{diag}

\DeclareMathOperator{\bmin}{\mathbf{min}}

\newtheorem{theorem}{Theorem}

\newtheorem{proposition}[theorem]{Proposition}

\numberwithin{equation}{section} \numberwithin{theorem}{section}

\begin{document}

\title{Pattern Recognition on Oriented Matroids:
Critical Committees and Distance Signals
}

\author{Andrey O. Matveev}
\email{andrey.o.matveev@gmail.com}


\begin{abstract}
If $\mathfrak{V}(\boldsymbol{R})$
is the vertex sequence of a symmetric cycle $\boldsymbol{R}$ in the tope graph of a simple acyclic oriented matroid $\mathcal{M}$ on a $t$-element ground set, then the set $\bmin\mathfrak{V}(\boldsymbol{R})$ of minimal elements in the subposet~$\mathfrak{V}(\boldsymbol{R})$ of the tope poset of $\mathcal{M}$, based at the positive tope, is a critical committee for\hfill $\mathcal{M}$\hfill that\hfill votes\hfill for\hfill the\hfill base\hfill tope.\hfill We\hfill consider\hfill the\hfill sequence
$\boldsymbol{z}_{\!\boldsymbol{R}}$ $:=(\rho(R):\ R\in\mathfrak{V}(\boldsymbol{R}))$
of poset ranks of the elements from the vertex sequence of $\boldsymbol{R}$ as a fragment of
a signal with period $2t$ and relate the number of members of the committee $\bmin\mathfrak{V}(\boldsymbol{R})$ to the magnitudes of~$\lfloor\frac{t}{2}\rfloor$ components, with odd indices, of the discrete Fourier transform of the distance vector~$\boldsymbol{z}_{\!\boldsymbol{R}}$.
\end{abstract}

\maketitle

\pagestyle{myheadings}

\markboth{A.O.~MATVEEV}{PATTERN RECOGNITION ON ORIENTED MATROIDS}

\thispagestyle{empty}

\tableofcontents

\section{Introduction}

If $\mathfrak{V}$ is the vertex sequence of some path in a connected graph, then the corresponding sequence $(d(w,v):\ v\in\mathfrak{V})$ of distances between its vertices and a distinguished vertex $w$ of the graph can be extended in a natural way to a periodic signal.
In this note, which is a companion to~\cite{AM1}, we consider
such distance signals associated with symmetric cycles in the tope graph of an oriented matroid. See~\cite{BLSWZ} on oriented matroids. When the vertex sequence of a symmetric cycle is regarded as a subposet of the tope poset and
the distinguished vertex of the tope graph is the
base tope of the tope poset, a few basic observations (see, e.g.,~\cite[Chapter~2]{F},\cite[Chapters~1$\div$6]{W}) concerning the discrete Fourier transform~(DFT) allow us to express the number of minimal elements of the vertex sequence of the symmetric cycle via the magnitudes of components of the~DFT of the distance signal. In the case of an acyclic oriented matroid on a $t$-element ground set, with the distinguished positive tope, we thus relate the number of members of a critical committee
to the magnitudes of $\lfloor\frac{t}{2}\rfloor$ components, with odd indices, of the DFT
of the distance signal.

\section{Critical Committees and Distance Signals}

Let
$\mathcal{A}:=(E_t,\mathcal{T})$
be a simple acyclic oriented matroid on the ground set~$E_t$ $:=\{1,\ldots,t\}$, with set of topes $\mathcal{T}$ whose components $-$ and $+$ are replaced by the real numbers $-1$ and $1$, respectively.
Let $\boldsymbol{R}:=(R^0,R^1,\ldots,R^{2t-1},R^0)$ be a symmetric cycle in the tope graph of $\mathcal{A}$, that is,
$R^{k+t}=-R^k$, $0\leq k$ $\leq t-1$. Denote by $\mathcal{T}(\mathcal{L}(\mathcal{A}),\mathrm{T}^{(+)})$ the tope poset of $\mathcal{A}$ based at the positive tope~$\mathrm{T}^{(+)}$. The set $\bmin\mathfrak{V}(\boldsymbol{R})$ of minimal elements of the subposet $\mathfrak{V}(\boldsymbol{R})$ $:=(R^0,R^1,\ldots,R^{2t-1})\subset\mathcal{T}(\mathcal{L}(\mathcal{A}),\mathrm{T}^{(+)})$
is a {\em critical tope committee for $\mathcal{A}$}: it is the inclusion-minimal subset $\mathcal{K}^{\ast}\subset\mathfrak{V}(\boldsymbol{R})$ such that
$\sum_{T\in\mathcal{K}^{\ast}}T=\mathrm{T}^{(+)}$.
Let $\rho(T)$ denote the poset rank of a tope $T\in\mathcal{T}(\mathcal{L}(\mathcal{A}),\mathrm{T}^{(+)})$.
The sequence
\begin{equation*}
\boldsymbol{z}_{\!\boldsymbol{R}}:=\bigl(z_{\!\boldsymbol{R}}(0):=\rho(R^0),z_{\!\boldsymbol{R}}(1):=\rho(R^1),\ldots,z_{\!\boldsymbol{R}}(2t-1):=\rho(R^{2t-1})\bigr)
\end{equation*}
determines the 
{\em distance signal\/}
$z_{\!\boldsymbol{R}}:\ \mathbb{Z}\to \{0\}\dot\cup E_t$
{\em of $\boldsymbol{R}$}, with period $2t$:
\begin{equation*}
z_{\!\boldsymbol{R}}(j+2t)=z_{\!\boldsymbol{R}}(j)\; ,\ \ \ j\in\mathbb{Z}\; .
\end{equation*}

Let
$\ell^2(\mathbb{Z}_{2t})$
denote the $2t$-dimensional complex coordinate space; the elements of $\ell^2(\mathbb{Z}_{2t})$ are supposed to be row vectors whose components are indexed from $0$ to $2t-1$.
We consider the
{\em distance vector\/} $\boldsymbol{z}_{\!\boldsymbol{R}}$ {\em of $\boldsymbol{R}$} as
an element of the space $\ell^2(\mathbb{Z}_{2t})$.

Let $\boldsymbol{\iota}$ denote the vector $(1,1,\ldots,1)\in\ell^2(\mathbb{Z}_{2t})$.
Since
$z_{\!\boldsymbol{R}}(k)+z_{\!\boldsymbol{R}}(k+t)=t$,
$0\leq k\leq t-1$,
we have
\begin{equation*}
\boldsymbol{z}_{\!\boldsymbol{R}}\cdot\boldsymbol{\iota}^{\!\!\top}=t^2\; .
\end{equation*}

If $R\in\mathfrak{V}(\boldsymbol{R})$ is a vertex of the symmetric cycle $\boldsymbol{R}$ then we denote
by~$\mathcal{N}(R)$ the
neighborhood of $R$ in the cycle $\boldsymbol{R}$;
$\mathbf{I}$ and $\mathbf{C}$ denote the $2t\times 2t$ identity matrix and basic circulant permutation matrix, respectively, with the rows and columns indexed from $0$ to $2t-1$.

$\bullet$ On the one hand,
$|\bmin\mathfrak{V}(\boldsymbol{R})|=t-\frac{1}{4}\sum_{\substack{\{R',R''\}:=\mathcal{N}(R):\\ R\in\mathfrak{V}(\boldsymbol{R})}}|\rho(R'')-\rho(R')|$, and we have
\begin{align*}
|\bmin\mathfrak{V}(\boldsymbol{R})|&=
t-\frac{1}{8}\sum_{\substack{\{R',R''\}:=\mathcal{N}(R):\\ R\in\mathfrak{V}(\boldsymbol{R})}}\bigl(\rho(R'')-\rho(R')\bigr)^2\\&=
t-\frac{1}{4}\sum_{\substack{\{R',R''\}:=\mathcal{N}(R):\\ R\in\mathfrak{V}(\boldsymbol{R})}}\bigl(\rho(R)^2-\rho(R')\rho(R'')\bigr)\; ;
\end{align*}
in other words, $|\bmin\mathfrak{V}(\boldsymbol{R})|=t-\frac{1}{4}\boldsymbol{z}_{\!\boldsymbol{R}}
\cdot(\mathbf{I}-\mathbf{C})(\mathbf{I}+\mathbf{C})\cdot
\boldsymbol{z}_{\!\boldsymbol{R}}{}^{\!\top}$, or
\begin{equation}
\label{eq:12}
|\bmin\mathfrak{V}(\boldsymbol{R})|=t-\frac{1}{8}\boldsymbol{z}_{\!\boldsymbol{R}}
\cdot(2\mathbf{I}-\mathbf{C}^{-2}-\mathbf{C}^2)\cdot
\boldsymbol{z}_{\!\boldsymbol{R}}{}^{\!\top}\; .
\end{equation}
The first row of the symmetric circulant matrix
$2\mathbf{I}-\mathbf{C}^{-2}-\mathbf{C}^2$ is the vector $\mathbf{b}:=(2,0,-1,0,\ldots,0,-1,0)$;
the components of the DFT of $\mathbf{b}$ are
\begin{align*}
\hat{\mathrm{b}}(k):=\sum_{n=0}^{2t-1} \mathrm{b}(n)\mathtt{e}^{-\pi\imath k n/t}=4\sin^2\tfrac{\pi k}{t}\; ,\ \ \ 0\leq k\leq 2t-1\; .
\end{align*}
Denote by $\mathbf{W}$ the $2t\times 2t$ Fourier matrix; its $(m,n)$-th entries are $\mathtt{e}^{-\pi\imath m n/t}$, $0\leq m,n\leq 2t-1$.
The DFT and inverse DFT of the distance vector $\boldsymbol{z}_{\!\boldsymbol{R}}$ are the vectors
$\hat{\boldsymbol{z}_{\!\boldsymbol{R}}}:=\boldsymbol{z}_{\!\boldsymbol{R}}\mathbf{W}$ and $\check{\boldsymbol{z}_{\!\boldsymbol{R}}}:=\boldsymbol{z}_{\!\boldsymbol{R}}\mathbf{W}^{-1}$, respectively;
thus,
\begin{equation*}
\begin{split}
\hat{z_{\!\boldsymbol{R}}}(k):&=\sum_{j=0}^{2t-1}z_{\!\boldsymbol{R}}(j)\mathtt{e}^{-\pi\imath k j/t}\\&=
\sum_{j=0}^{t-1}\bigl(z_{\!\boldsymbol{R}}(j)\mathtt{e}^{-\pi\imath k j/t} +
(t-z_{\!\boldsymbol{R}}(j))\mathtt{e}^{-\pi\imath k (j+t)/t}\bigr)\; ,\ \ \ 0\leq k\leq 2t-1\; ,
\end{split}
\end{equation*}
that is,
\begin{equation*}
\hat{z_{\!\boldsymbol{R}}}(k)=
\begin{cases}
t^2\; ,& \text{$k=0$}\; ,\\
0\; ,& \text{$k$ even, $k\neq 0$}\; ,\\
2\bigl(-t(1-\mathtt{e}^{-\pi\imath k/t})^{-1}+\sum_{j=0}^{t-1}z_{\!\boldsymbol{R}}(j)\mathtt{e}^{-\pi\imath k j/t}\bigr)\; ,& \text{$k$ odd}\; ;
\end{cases}
\end{equation*}
in particular, if $t$ is odd then
$\hat{z_{\!\boldsymbol{R}}}(t)=-t+2\sum_{j=0}^{t-1}(-1)^j z_{\!\boldsymbol{R}}(j)$.

We have
\begin{equation*}
2\mathbf{I}-\mathbf{C}^{-2}-\mathbf{C}^2=
\mathbf{W}^{-1}\cdot 4\diag(0,\sin^2\tfrac{\pi}{t},\ldots,\sin^2\tfrac{\pi k}{t},\ldots,\sin^2\tfrac{\pi (2t-1)}{t})\cdot\mathbf{W}\; ,
\end{equation*}
and Eq.~(\ref{eq:12}) implies that
\begin{equation}
\label{eq:14}
|\bmin\mathfrak{V}(\boldsymbol{R})|=
t-\tfrac{1}{2}\check{\boldsymbol{z}_{\!\boldsymbol{R}}}
\cdot \diag(0,\sin^2\tfrac{\pi}{t},\ldots,\sin^2\tfrac{\pi k}{t},\ldots,\sin^2\tfrac{\pi (2t-1)}{t})\cdot\hat{\boldsymbol{z}_{\!\boldsymbol{R}}}^{\top}\; .
\end{equation}
Denote by $\overline{\hat{\boldsymbol{z}_{\!\boldsymbol{R}}}}$ the vector composed of the complex conjugates of components of $\hat{\boldsymbol{z}_{\!\boldsymbol{R}}}$.
Since $\check{\boldsymbol{z}_{\!\boldsymbol{R}}}=\tfrac{1}{2t}\overline{\hat{\boldsymbol{z}_{\!\boldsymbol{R}}}}$, it follows from  Eq.~(\ref{eq:14}) that
\begin{equation*}
|\bmin\mathfrak{V}(\boldsymbol{R})|=
t-\tfrac{1}{4t}\overline{\hat{\boldsymbol{z}_{\!\boldsymbol{R}}}}
\cdot \diag(0,\sin^2\tfrac{\pi}{t},\ldots,\sin^2\tfrac{\pi k}{t},\ldots,\sin^2\tfrac{\pi (2t-1)}{t})\cdot\hat{\boldsymbol{z}_{\!\boldsymbol{R}}}^{\top}\; .
\end{equation*}
Using Plancherel's formula, we reformulate this observation:
\begin{align}
\label{eq:18}
|\bmin\mathfrak{V}(\boldsymbol{R})|&=
t-\frac{1}{4t}\sum_{k=0}^{2t-1} |\hat{z_{\!\boldsymbol{R}}}(k)|^2\cdot\sin^2\tfrac{\pi k}{t}
\\ \nonumber
&=t-\frac{1}{2}\|\boldsymbol{z}_{\!\boldsymbol{R}}\|^2+
\frac{1}{4t}\sum_{k=0}^{2t-1} |\hat{z_{\!\boldsymbol{R}}}(k)|^2\cdot\cos^2\tfrac{\pi k}{t}\; .
\end{align}
Since $\hat{z_{\!\boldsymbol{R}}}(k)=\overline{\hat{z_{\!\boldsymbol{R}}}(2t-k)}$ and the magnitudes $|\hat{z_{\!\boldsymbol{R}}}(k)|$ do not depend on circular translation:
$|(\boldsymbol{z}_{\!\boldsymbol{R}}\mathbf{C}^j)\hat{}\;(k)|=|\hat{z_{\!\boldsymbol{R}}}(k)|$, $1\leq k\leq 2t-1$, for any integer~$j$,
we restate Eq.~(\ref{eq:18}) in the following way:
\begin{proposition}
\label{prop:2}
\begin{itemize}
\item[(\rm{i})]
For any distance vector  $\boldsymbol{z}_{\!\boldsymbol{R}}$ of the symmetric cycle $\boldsymbol{R}$, we have
\begin{equation}
\label{eq:19}
|\bmin\mathfrak{V}(\boldsymbol{R})|=
t-\frac{1}{2t}\sum_{\substack{1\leq k\leq t-1,\\ \text{\rm $k$ odd}}} |\hat{z_{\!\boldsymbol{R}}}(k)|^2\cdot\sin^2\tfrac{\pi k}{t}\; .
\end{equation}
\item[(\rm{ii})]
Let $\mathcal{M}:=(E_t,\mathcal{T})$ be a simple oriented matroid {\rm(}it has no loops, parallel or {\sl antiparallel\/}  elements{\rm)}, and $\mathfrak{V}(\boldsymbol{R}):=(R^0,R^1,\ldots,R^{2t-1})$
the vertex sequence of a
symmetric cycle $\boldsymbol{R}$ in the tope graph of $\mathcal{M}$.
Given a tope $T$ of $\mathcal{M}$, denote by $\boldsymbol{z}_{T,\boldsymbol{R}}$ the {\em distance vector of
$\boldsymbol{R}$ with respect to the tope $T$}, that is, the sequence $(d(T,R):\ R\in\mathfrak{V}(\boldsymbol{R}))$
of
graph
distances between $T$ and the vertices of $\boldsymbol{R}$. For the inclusion-minimal subset $\boldsymbol{Q}(T,\boldsymbol{R})\subset\mathfrak{V}(\boldsymbol{R})$ such that $T=\sum_{Q\in\boldsymbol{Q}(T,\boldsymbol{R})}Q$, we have
\begin{equation}
\label{eq:20}
|\boldsymbol{Q}(T,\boldsymbol{R})|=
t-\frac{1}{2t}\sum_{\substack{1\leq k\leq t-1,\\ \text{\rm $k$ odd}}} |\hat{z_{T}}{}_{\!,\boldsymbol{R}}(k)|^2\cdot\sin^2\tfrac{\pi k}{t}\; .
\end{equation}
\end{itemize}
\end{proposition}
{\small
The vertex sequence $\mathfrak{V}(\boldsymbol{R})$ is a maximal positive basis of $\mathbb{R}^t$; the existence of the set $\boldsymbol{Q}(T,\boldsymbol{R})$,
of odd cardinality, of linearly independent elements of $\mathbb{R}^t$, mentioned in Proposition~\ref{prop:2}(ii), is guaranteed by~\cite[Corollary~2.2]{AM1}.

Following~\cite[Lecture~6]{St}, define the {\em distance enumerator}, with respect to a tope~$B\in\mathcal{T}$, of the set $\mathcal{T}$, as the polynomial $D_{B,\mathcal{T}}(\mathrm{x}):=\sum_{T\in\mathcal{T}}\mathrm{x}^{d(B,T)}$. In an analogous manner, define the distance enumerator of the vertex set of the symmetric cycle $\boldsymbol{R}$ as the summand $D_{B,\mathfrak{V}(\boldsymbol{R})}(\mathrm{x}):=\sum_{T\in\mathfrak{V}(\boldsymbol{R})}\mathrm{x}^{d(B,T)}$ of $D_{B,\mathcal{T}}(\mathrm{x})$.

Note that for any topes $T'$ and $T''$ of $\mathcal{M}$ the graph distance between them is
$d(T',T'')$ $=t-\tfrac{1}{4}\|T''+T'\|^2 =\tfrac{1}{4}\|T''-T'\|^2=\tfrac{1}{2}\bigl(t-\langle T'',T'\rangle\bigr)$. Thus, if $\mathcal{T}^{\bullet}$ is a halfspace of $\mathcal{M}$, and if $\mathfrak{W}\subset\mathfrak{V}(\boldsymbol{R})$ is the vertex sequence of a $(t-1)$-path in the cycle $\boldsymbol{R}$, then
$D_{B,\mathcal{T}}(\mathrm{x})=\mathrm{x}^{t/2}\sum_{T\in\mathcal{T}^{\bullet}}
(\mathrm{x}^{-\langle B,T\rangle/2} + \mathrm{x}^{\langle B,T\rangle/2})$
and
$D_{B,\mathfrak{V}(\boldsymbol{R})}(\mathrm{x})$ $=\mathrm{x}^{t/2}\sum_{T\in\mathfrak{W}}
(\mathrm{x}^{-\langle B,T\rangle/2} + \mathrm{x}^{\langle B,T\rangle/2})$.

Associate with a tope $T\in\mathcal{T}$ the $2t$-dimensional row vector $\boldsymbol{q}(T,\boldsymbol{R})$ defined by~$q_j(T,\boldsymbol{R}):=1$ if $R^j\in\boldsymbol{Q}(T,\boldsymbol{R})$, and $q_j(T,\boldsymbol{R}):=0$ otherwise. Let $\mathbf{G}(\boldsymbol{R})$ denote the Gram matrix of the sequence $\mathfrak{V}(\boldsymbol{R})$. For topes $T'$ and $T''$ of $\mathcal{M}$ we have
$d(T',T'')=\tfrac{1}{2}\bigl(t-\boldsymbol{q}(T'',\boldsymbol{R})\mathbf{G}(\boldsymbol{R})\boldsymbol{q}(T',\boldsymbol{R})^{\top}\bigr)$.
}

$\bullet$ On the other hand,
\begin{equation*}
\begin{split}
|\bmin\mathfrak{V}(\boldsymbol{R})|&=
\frac{1}{8}\sum_{\substack{\{R',R''\}:=\mathcal{N}(R):\\ R\in\mathfrak{V}(\boldsymbol{R})}}\left(\rho(R')+\rho(R'')-2\rho(R)\right)^2\\&=
\frac{1}{4}\sum_{\substack{\{R',R''\}:=\mathcal{N}(R):\\ R\in\mathfrak{V}(\boldsymbol{R})}}\left(
3\rho(R)^2-2\rho(R)(\rho(R')+\rho(R''))+\rho(R')\rho(R'')
\right)\; ,
\end{split}
\end{equation*}
that is, $|\bmin\mathfrak{V}(\boldsymbol{R})|=\frac{1}{4}\boldsymbol{z}_{\!\boldsymbol{R}}
\cdot(\mathbf{I}-\mathbf{C})(3\mathbf{I}-\mathbf{C})\cdot
\boldsymbol{z}_{\!\boldsymbol{R}}{}^{\!\top}$, or
\begin{equation}
\label{eq:15}
|\bmin\mathfrak{V}(\boldsymbol{R})|=
\frac{1}{8}\boldsymbol{z}_{\!\boldsymbol{R}}
\cdot(6\mathbf{I}-4\mathbf{C}^{-1}-4\mathbf{C}+\mathbf{C}^{-2}+\mathbf{C}^2)\cdot
\boldsymbol{z}_{\!\boldsymbol{R}}{}^{\!\top}
\; .
\end{equation}

Eqs.~(\ref{eq:12}) and~(\ref{eq:15}) imply that
$|\bmin\mathfrak{V}(\boldsymbol{R})|=\frac{t}{2}+\frac{1}{4}\boldsymbol{z}_{\!\boldsymbol{R}}
\cdot(\mathbf{I}-\mathbf{C})^2\cdot
\boldsymbol{z}_{\!\boldsymbol{R}}{}^{\!\top}$ $=\frac{3t}{4}-\frac{1}{4}\boldsymbol{z}_{\!\boldsymbol{R}}
\cdot(\mathbf{I}-\mathbf{C})\mathbf{C}\cdot
\boldsymbol{z}_{\!\boldsymbol{R}}{}^{\!\top}$; in other words, we have
\begin{align}
\label{eq:16}
|\bmin\mathfrak{V}(\boldsymbol{R})|&=
\frac{t}{2}+\frac{1}{8}\boldsymbol{z}_{\!\boldsymbol{R}}
\cdot(2\mathbf{I}-2\mathbf{C}^{-1}-2\mathbf{C}+\mathbf{C}^{-2}+\mathbf{C}^2)\cdot
\boldsymbol{z}_{\!\boldsymbol{R}}{}^{\!\top}
\intertext{and}
\label{eq:17}
|\bmin\mathfrak{V}(\boldsymbol{R})|&=\frac{3t}{4}-\frac{1}{8}\boldsymbol{z}_{\!\boldsymbol{R}}
\cdot(\mathbf{C}^{-1}+\mathbf{C}-\mathbf{C}^{-2}-\mathbf{C}^2)\cdot
\boldsymbol{z}_{\!\boldsymbol{R}}{}^{\!\top}\; ;
\end{align}
note also that
\begin{equation*}
\boldsymbol{z}_{\!\boldsymbol{R}}
\cdot(\mathbf{I}-\mathbf{C})\cdot
\boldsymbol{z}_{\!\boldsymbol{R}}{}^{\!\top}=t\; .
\end{equation*}
We 
derive from Eqs.~(\ref{eq:15}), (\ref{eq:16}) and~(\ref{eq:17}) the relations
\begin{align*}
|\bmin\mathfrak{V}(\boldsymbol{R})|&=\frac{1}{4t}\sum_{k=0}^{2t-1} |\hat{z_{\!\boldsymbol{R}}}(k)|^2\cdot
(\cos^2\tfrac{\pi k}{t}-2\cos\tfrac{\pi k}{t}+1)\; ,\\
|\bmin\mathfrak{V}(\boldsymbol{R})|&=\frac{t}{2}+
\frac{1}{4t}\sum_{k=0}^{2t-1} |\hat{z_{\!\boldsymbol{R}}}(k)|^2\cdot
(\cos^2\tfrac{\pi k}{t}-\cos\tfrac{\pi k}{t})
\intertext{and}
|\bmin\mathfrak{V}(\boldsymbol{R})|&=\frac{3t}{4}+
\frac{1}{8t}\sum_{k=0}^{2t-1} |\hat{z_{\!\boldsymbol{R}}}(k)|^2\cdot
(2\cos^2\tfrac{\pi k}{t}-\cos\tfrac{\pi k}{t}-1)
\end{align*}
that are equivalent to Eq.~(\ref{eq:18}).

$\bullet$ Let $\boldsymbol{R}'$ and $\boldsymbol{R}''$ be two symmetric cycles in the tope graph of the oriented matroid\hfill $\mathcal{A}$,\hfill and\hfill $\boldsymbol{z}_{\!\boldsymbol{R}'}$\hfill and\hfill $\boldsymbol{z}_{\!\boldsymbol{R}''}$\hfill their\hfill distance\hfill vectors.\hfill Consider\hfill the\newline
vectors $\boldsymbol{e}:=\boldsymbol{z}_{\!\boldsymbol{R}''}-\boldsymbol{z}_{\!\boldsymbol{R}'}$ and $\boldsymbol{m}:=\boldsymbol{z}_{\!\boldsymbol{R}''}+\boldsymbol{z}_{\!\boldsymbol{R}'}$.
Note that the components of their DFTs are:
\begin{align*}
\hat{e}(k)&=
\begin{cases}
0\; ,& \text{$k$ even}\; ,\\
2\sum_{j=0}^{t-1}e(j)\mathtt{e}^{-\pi\imath k j/t}\; ,& \text{$k$ odd}\; ,
\end{cases}
\intertext{and}
\hat{m}(k)&=
\begin{cases}
2t^2\; ,& \text{$k=0$}\; ,\\
0\; ,& \text{$k$ even, $k\neq 0$}\; ,\\
2\bigl(-2t(1-\mathtt{e}^{-\pi\imath k/t})^{-1}+\sum_{j=0}^{t-1}m(j)\mathtt{e}^{-\pi\imath k j/t}\bigr)\; ,& \text{$k$ odd}\; ;
\end{cases}
\end{align*}
thus, if $t$ is odd then $\hat{e}(t)=2\sum_{j=0}^{t-1}(-1)^j e(j)$
and $\hat{m}(t)=2\bigl(-t+\sum_{j=0}^{t-1}(-1)^j m(j)\bigr)$.

It follows from Proposition~\ref{prop:2}(i) that
\begin{equation*}
|\bmin\mathfrak{V}(\boldsymbol{R}')|+|\bmin\mathfrak{V}(\boldsymbol{R}'')|=
2t-\frac{1}{4t}\sum_{\substack{1\leq k\leq t-1,\\ \text{\rm $k$ odd}}}\bigl(|\hat{e}(k)|^2+
|\hat{m}(k)|^2\bigr)\cdot\sin^2\tfrac{\pi k}{t}\; .
\end{equation*}

\end{document}